\documentclass[11pt]{article}

\usepackage[T1]{fontenc}
\usepackage[utf8]{inputenc}
\usepackage[american]{babel}
\usepackage[margin=1in]{geometry}

\usepackage{amsmath}
\usepackage{amssymb}
\usepackage{amsthm}
\usepackage{dsfont}
\usepackage{color}
\usepackage{url}
\usepackage{xcolor}
\usepackage{enumitem}
\usepackage{booktabs}
\definecolor{darkgreen}{rgb}{0.13, 0.55, 0.13}
\definecolor{C0}{HTML}{3182ce}
\usepackage[breaklinks=true,colorlinks,bookmarks=false,citecolor=C0,linkcolor=C0,urlcolor=C0]{hyperref}

\renewcommand{\Pr}{\operatorname{pr}}

\usepackage{natbib}

\theoremstyle{plain}
\newtheorem{theorem}{Theorem}
\newtheorem{lemma}{Lemma}
\newtheorem{proposition}{Proposition}
\newtheorem{corollary}{Corollary}
\newtheorem{factbase}{Fact}
\theoremstyle{definition}
\newtheorem{definition}{Definition}
\newtheorem{example}{Example}
\theoremstyle{remark}
\newtheorem{remark}{Remark}

\usepackage{xparse}
\NewDocumentEnvironment{thm}{o}
  {\IfNoValueTF{#1}{\begin{theorem}}{\begin{theorem}[#1]}}
  {\end{theorem}}
\NewDocumentEnvironment{lem}{o}
  {\IfNoValueTF{#1}{\begin{lemma}}{\begin{lemma}[#1]}}
  {\end{lemma}}
\NewDocumentEnvironment{prop}{o}
  {\IfNoValueTF{#1}{\begin{proposition}}{\begin{proposition}[#1]}}
  {\end{proposition}}
\NewDocumentEnvironment{cor}{o}
  {\IfNoValueTF{#1}{\begin{corollary}}{\begin{corollary}[#1]}}
  {\end{corollary}}
\NewDocumentEnvironment{rem}{o}
  {\IfNoValueTF{#1}{\begin{remark}}{\begin{remark}[#1]}}
  {\end{remark}}
\NewDocumentEnvironment{defn}{o}
  {\IfNoValueTF{#1}{\begin{definition}}{\begin{definition}[#1]}}
  {\end{definition}}
\NewDocumentEnvironment{fact}{o}
  {\IfNoValueTF{#1}{\begin{factbase}}{\begin{factbase}[#1]}}
  {\end{factbase}}

\newenvironment{keywords}{%
  \par\smallskip\noindent\textbf{Keywords:}\ }{\par\smallskip}

\begin{document}
\raggedbottom

\title{Consistency of variational approximations\\ under bounded Kullback--Leibler divergence}

\author{%
Hien Duy Nguyen\\
Department of Mathematics and Physical Science,\\ La Trobe University, Melbourne, Australia\\
Institute of Mathematics for Industry, Kyushu University, Fukuoka, Japan\\
\texttt{hien@imi.kyushu-u.ac.jp}
\and
Jacob Westerhout\\
School of Mathematics and Physics, University of Queensland, Brisbane, Australia\\
\texttt{j.westerhout@uq.edu.au}
\and
Thomas Guilmeau and Julyan Arbel\\
Univ. Grenoble Alpes, Inria, CNRS, Grenoble INP, LJK, 38000 Grenoble, France\\
\texttt{thomas.guilmeau@inria.fr}, \texttt{julyan.arbel@inria.fr}}

\date{}

\maketitle

\begin{abstract}
Variational methods are widely used to approximate posterior distributions in Bayesian inference when exact computation is infeasible. We study when such approximations inherit posterior consistency. Our first result shows that, on a general metric space, a uniform bound on the Kullback--Leibler divergence from the approximating measures to a tight sequence of target measures forces the approximating sequence to be tight. It follows that if the target posteriors converge weakly to a Dirac mass at the true parameter, then any variational sequence with bounded Kullback--Leibler divergence to the targets is also consistent. We also give simple logarithmic-moment conditions that verify this boundedness condition, and illustrate them for smooth generalised posterior distributions.
\end{abstract}


\section{Introduction}

Variational methods are widely used to approximate posterior probability measures arising in Bayesian and generalised Bayesian inference, particularly when exact computation is infeasible \citep{blei2017variational, zhang2019}.
A posterior sequence is consistent when the posterior distributions $\mu_n$ converge to a Dirac mass at the true parameter $\theta_0$ as the number of observations $n$ grows to infinity.
A fundamental theoretical question is as follows: under what conditions do variational approximations $\nu_n$ inherit this convergence?

This question is the focus of a growing literature. However, most existing works impose restrictions on one or more of the following aspects: (i) the approximating family, (ii) the statistical model or posterior, and (iii) the variational optimisation procedure. \citet{cherief2018consistency-mixtures} study mean-field variational Bayes for mixture estimation and model selection, while \citet{RaySzabo2022VBHDLR} study mean-field spike-and-slab approximations for high-dimensional sparse linear regression. \citet{cherief-abdellatif2020convergence} treats sparse deep learning and model selection over neural-network architectures. \citet{alquier2020concentration} and \citet{yang2020alpha} work with fractional or $\alpha$-variational posteriors; the latter extends the framework to latent-variable models. \citet{whang2019} obtain asymptotic variational Bernstein--von Mises results for parametric latent-variable models, and \citet{zhang2020convergence} give convergence rates under prior-mass and testing conditions, including conditions tailored to mean-field approximations. These contributions provide statistical guarantees in important settings, but their assumptions typically encode some combination of model structure, variational-family structure, exact variational optimisation, or fractional posterior form. As a result, the basic mechanism by which posterior consistency is preserved under variational approximation remains unclear.

A typical variational inference construction produces a sequence of approximations by (approximately) minimising the Kullback--Leibler (KL) divergence with respect to the posterior $\mu_n$ over a variational class. In this paper, we show that uniform control of this KL divergence induces tightness of variational sequences, leading to a general consistency principle in metric spaces. The resulting argument, building on a technical device introduced by \citet{RaySzabo2022VBHDLR}, is concise and relies only on fundamental properties of the KL divergence and weak convergence of probability measures. In particular, we impose no structural assumption on the approximating family, the possibly tempered posteriors only need to converge to a Dirac mass, and variational approximations that only approximately minimise the KL divergence are allowed.

To illustrate the usefulness of this principle, we derive simple sufficient conditions ensuring boundedness of the KL divergence based on moment bounds, and apply these results to smooth generalised posterior distributions. Thus, consistency of variational approximations can be reduced to verifying an interpretable moment condition, rather than a detailed analysis of the variational optimisation problem itself. These results hold for a large class of variational families, including in particular location-scale families, such as Gaussian variational approximations.

In Section~\ref{sec:metric-space}, we first establish a general tightness property for KL level sets in metric spaces. We show in Theorem~\ref{thm:metric:asympt-tight} that uniform KL control relative to a tight sequence of target measures induces tightness of the corresponding approximating sequence.
Our main contribution, Theorem~\ref{thm:metric:delta}, establishes consistency for variational sequences and follows directly from Theorem~\ref{thm:metric:asympt-tight}, under a sufficient boundedness condition on the KL divergence between $\nu_n$ and $\mu_n$. We then use three examples to show that natural weakenings of this boundedness condition do not lead to consistency and that this condition is not necessary.
Section~\ref{sec:log-moments} addresses the complementary question of verifying the KL boundedness condition in practice. We show in Proposition~\ref{prop:log-mom-oneQ} that boundedness of the optimal variational KL divergence can be reduced to a simple moment condition under a single reference distribution. We then provide a sufficient condition in Proposition~\ref{prop:LAN-like} based on a quadratic lower envelope, while Proposition~\ref{prop:quadratic-envelope} shows how to handle this envelope condition. We illustrate the approach on generalised posterior distributions in Theorem~\ref{thm:gp-envelope}. Finally, Section~\ref{sec:GLM} applies these results in the context of generalised linear models.\\

\noindent\textbf{Notations.} Let $(\mathcal X,d)$ be a metric space and let $\mathcal P(\mathcal X)$ denote the set of Borel probability measures on $\mathcal X$.
We write weak convergence as $\nu_n\rightsquigarrow\nu$.
For $\mu,\nu\in\mathcal P(\mathcal X)$, the KL divergence is
$\mathrm{KL}(\nu\Vert\mu)=\int_{\mathcal X}\log({\mathrm{d}\nu}/{\mathrm{d}\mu})\,\mathrm{d}\nu$,
with $\mathrm{KL}(\nu\Vert\mu)=\infty$ if $\nu\not\ll\mu$.
We say that $\mu$ and $\nu$ are equivalent if $\mu\ll\nu$ and $\nu\ll\mu$.
For $\delta>0$ and $A\subseteq \mathcal X$, we denote the $\delta$-expansion of $A$ by $A^\delta = \{x\in\mathcal X: d(x,A) < \delta\}$, and denote the $\delta$-shrinkage of $A$ by $A_\delta = \{x\in A: B(x,\delta)\subseteq A\}$, where $B(x,\delta)=\{y\in\mathcal X:d(y,x)<\delta\}$ is the open ball of radius $\delta$ centred at $x$. We write $\overline{B}(x,\delta)$ for the corresponding closed ball. Note that $(A^c)_\delta = (A^\delta)^c$. The indicator function of a set $A$ is denoted by $\mathds{1}_A$, and $x_+=\max(x,0)$ for $x\in\mathbb R$.

\section{Consistency of variational approximations in metric spaces through tightness}
\label{sec:metric-space}

This section collects the main results on consistency of variational sequences in metric spaces. The idea of the following technical device is due to \citet{RaySzabo2022VBHDLR}.

\begin{lem}\label{lem:j:1}
For probability measures $\mu,\nu$ on a common measurable space $(\mathcal{X},\mathfrak X)$, for every $A\in\mathfrak X$ and every $\delta>0$,
\[
\nu(A) \le \frac{1}{\delta}\left\{\mathrm{KL}(\nu\Vert\mu) + \mu(A) e^{\delta}\right\}.
\]
\end{lem}

\begin{proof}
If $\nu\not\ll \mu$ then $\mathrm{KL}(\nu\Vert\mu)=\infty$ and the bound is trivial.
Otherwise, the Donsker--Varadhan variational bound \cite[Thm. 4.6]{polyanskiy_wu_2025_information_theory} implies that for any measurable $f$,
\[
\int f\,\mathrm{d}\nu \le \mathrm{KL}(\nu\Vert\mu) + \ln\left(\int e^{f}\,\mathrm{d}\mu\right).
\]
Taking $f=\delta\,\mathds{1}_A$ and using $\ln(1+x)\le x$ yields the result.
\end{proof}

\subsection{A KL level-set tightness statement}

\begin{definition}
    Let $(\mathcal X,d)$ be a metric space and for every $n\in \mathbb N$ let $S_n$ be a set of Borel probability measures.
    We say that $(S_n)$ is asymptotically tight if for every $\eta>0$ there exists a compact $K\subseteq\mathcal X$ such that for every sequence $(\nu_n)$ with $\nu_n\in S_n$ and every $\delta>0$,
    \[
    \liminf_{n\to\infty}\nu_n(K^\delta)\ge 1-\eta.
    \]
    When every $S_n$ is a singleton, this reduces to the standard definition of an asymptotically tight sequence of probability measures \cite[Def 1.3.7]{vaart2023empirical}.
\end{definition}

\begin{thm}\label{thm:metric:asympt-tight}
Let $(\mathcal X,d)$ be a metric space and let $(\mu_n)$ be a sequence of asymptotically tight Borel probability measures on $\mathcal X$.
For any sequence $(\epsilon_n)\subseteq [0,\infty]$ with $\limsup_n \epsilon_n<\infty$, let the sets
\[
S_n=\{\nu\in\mathcal P(\mathcal X):\ \mathrm{KL}(\nu\Vert\mu_n)\le \epsilon_n\}.
\]
Then the sequence $(S_n)$ is asymptotically tight.
\end{thm}

\begin{proof}
Let $M=\limsup_n\epsilon_n<\infty$ and fix $\eta>0$.
Choose $\zeta>0$ such that $(M+1)/\ln(1/\zeta)\le\eta$.
By asymptotic tightness of $(\mu_n)$ there exists a compact $K\subseteq\mathcal X$ such that, for every $a>0$,
\[
\limsup_{n\to\infty} \mu_n\left((K^c)_a\right)\le\zeta.
\]
Let $(\nu_n)$ be any sequence with $\nu_n\in S_n$. For fixed $a>0$, apply Lemma~\ref{lem:j:1} with $A=(K^c)_a$ and with the scalar parameter $u=\ln(1/\zeta)$. Since $\mathrm{KL}(\nu_n\Vert\mu_n)\le\epsilon_n$,
\[
\limsup_{n\to\infty}\nu_n((K^c)_a)\le \frac{M+\zeta e^u}{u}=\frac{M+1}{\ln(1/\zeta)}\le\eta.
\]
Using $(K^c)_a=(K^a)^c$, this is equivalent to $\liminf_n\nu_n(K^a)\ge1-\eta$ for every $a>0$.
\end{proof}

\subsection{A consistency consequence}

\begin{thm}\label{thm:metric:delta}
Let $(\mathcal X,d)$ be a metric space and let $\mu_n\rightsquigarrow\delta_{\theta_0}$.
For any $\nu_n\in\mathcal P(\mathcal X)$ such that $\limsup_n \mathrm{KL}(\nu_n\Vert\mu_n)< \infty$, we have $\nu_n\rightsquigarrow\delta_{\theta_0}$.
\end{thm}

\begin{proof}
Since $\delta_{\theta_0}$ is tight, $(\mu_n)$ is asymptotically tight \cite[Lemma 1.3.8]{vaart2023empirical}.
Applying Theorem~\ref{thm:metric:asympt-tight} with $\epsilon_n=\mathrm{KL}(\nu_n\Vert\mu_n)$ shows that $(\nu_n)$ is asymptotically tight. Hence, by Prohorov's theorem in metric spaces \cite[Thm.~1.3.9]{vaart2023empirical}, every subsequence of $(\nu_n)$ has a further weakly convergent subsequence.
Let $\nu_{n_k}\rightsquigarrow\nu$ be such a subsequential limit. Lower semicontinuity of KL under weak convergence \cite[Thm.~4.9]{polyanskiy_wu_2025_information_theory} gives
\[
\mathrm{KL}(\nu\Vert\delta_{\theta_0}) \le \liminf_k \mathrm{KL}(\nu_{n_k}\Vert\mu_{n_k}) < \infty.
\]
Thus $\nu\ll\delta_{\theta_0}$, and hence $\nu=\delta_{\theta_0}$. Every subsequence of $(\nu_n)$ therefore has a further subsequence converging weakly to $\delta_{\theta_0}$, which implies $\nu_n\rightsquigarrow\delta_{\theta_0}$.
\end{proof}

Note that many natural weakenings of the condition $\limsup_{n}\mathrm{KL}(\nu_{n}\Vert\mu_{n})<\infty$
do not imply $\nu_{n}\rightsquigarrow\delta_{\theta_{0}}$.
Example~\ref{ex:1-better} shows that finiteness of  $\mathrm{KL}(\nu_{n}\Vert\mu_{n})$ for each
$n$ is not even sufficient for $\nu_n$ to converge.
The same example also shows that these failures can occur with $\mathrm{KL}(\nu_{n}\Vert\mu_{n})$ diverging arbitrarily slowly.
Note also that boundedness of the KL divergence is not necessary.
Example~\ref{ex:2}
satisfies $\mu_{n},\nu_{n}\rightsquigarrow\delta_{0}$ with $\mu_n$ equivalent to $\nu_n$, but $\mathrm{KL}(\nu_{n}\Vert\mu_{n})=\infty$
for every $n$.
Finally, Example~\ref{ex:3} shows that even if $\mu_{n},\nu_{n}\rightsquigarrow\delta_{0}$ and
$\mathrm{KL}(\nu_{n}\Vert\mu_{n})<\infty$ for every $n$, one may
still have $\limsup_{n\to\infty}\mathrm{KL}(\nu_{n}\Vert\mu_{n})=\infty$.

\begin{example}\label{ex:1-better}
Take any positive sequence $\epsilon_{n}\to\infty$, any sequence $(a_n)\subseteq \mathbb R$ and let $\mu_{n}$
and $\nu_{n}$ be the measures of laws ${\rm N}(0,(2\epsilon_{n})^{-1})$ and
${\rm N}(a_n,(2\epsilon_{n})^{-1})$, respectively. Then $\mu_{n}$ and $\nu_{n}$ are equivalent, $\mu_{n}\rightsquigarrow\delta_{0}$,
and $\mathrm{KL}(\nu_{n}\Vert\mu_{n})=\epsilon_{n}a_n^2<\infty$ for every $n$.
However, $\nu_{n}$ converges if and only if $a_n$ converges to, say, $a$, in which case $\nu_n \rightsquigarrow \delta_a$.

\end{example}

\begin{example} \label{ex:2}
Let $\phi_{\sigma}$ and $\psi$ denote, respectively, the densities of ${\rm N}(0,\sigma^{2})$ and the standard Cauchy distribution. Define
\[
q_{n}(x)=\left(1-\frac{1}{n}\right)\phi_{1/n}(x)+\frac{1}{n}\psi(x),\  p_{n}(x)=Z_{n}^{-1}e^{-|x|}q_{n}(x),\  Z_{n}=\int_{\mathbb{R}}e^{-|x|}q_{n}(x)\,dx\in(0,1),
\]
and let $\nu_{n}$ and $\mu_{n}$ have densities $q_{n}$ and $p_{n}$,
respectively. Then $\mu_{n}$ and $\nu_{n}$ are equivalent, $\nu_{n}\rightsquigarrow\delta_{0}$,
and, since $p_{n}\propto e^{-|x|}q_{n}$ with $e^{-|x|}$ bounded
and continuous, also $\mu_{n}\rightsquigarrow\delta_{0}$. Moreover,
\[
\mathrm{KL}(\nu_{n}\Vert\mu_{n})=\int\ln\!\left({q_{n}}/{p_{n}}\right)\,d\nu_{n}=\ln Z_{n}+\int|x|\,d\nu_{n}(x)=\infty,
\]
because $\nu_{n}$ has a Cauchy component of weight $1/n$.
\end{example}
\begin{example}\label{ex:3}
Let $\mu_{n}$ and $\nu_{n}$ be the laws of ${\rm N}(0,n^{-2})$
and ${\rm N}(0,n^{-1})$, respectively. Then $\mu_{n},\nu_{n}\rightsquigarrow\delta_{0}$,
whereas
\[
\mathrm{KL}(\nu_{n}\Vert\mu_{n})=\frac{1}{2}\left(n-1-\ln n\right)\to\infty.
\]
\end{example}

\subsection{Relation to the consistency of variational sequences}

The following theorem, which is a direct corollary of Theorem~\ref{thm:metric:delta}, states the sufficient Assumption (A1) which links the approximating families $\mathcal{F}_n$ and the target distributions $\mu_n$.

\begin{thm}
    Let $(\mathcal X,d)$ be a metric space, $\mu_n$ a sequence of Borel probability measures on $\mathcal X$ converging weakly to $\delta_{\theta_0}$ with $\theta_0 \in \mathcal X$.
    Let $\mathcal F_n$ denote a sequence of sets of Borel probability measures on $\mathcal X$ and
    assume
    \begin{itemize}[leftmargin=1.2cm]
        \item[$\mathrm{(A1)}$] $m_{n}=\inf_{\nu\in\mathcal{F}_{n}}\mathrm{KL}(\nu\Vert\mu_{n})$ satisfies $\limsup_n m_n < \infty$.
    \end{itemize}
    Then for any sequence $\nu_n \in \mathcal F_n$ with $\mathrm{KL}(\nu_n\Vert \mu_n) \leq m_n + \epsilon_n$ where $\limsup_n \epsilon_n<\infty$, $\nu_n$ converges weakly to $\delta_{\theta_0}$: $\nu_n\rightsquigarrow \delta_{\theta_0}$.
\end{thm}
\begin{proof}
    This is immediate by Theorem \ref{thm:metric:delta}.
\end{proof}

Conditions under which Assumption (A1) is satisfied are discussed in Section~\ref{sec:log-moments}. The assumption that $\limsup_n \epsilon_n < \infty$ depends on how $\nu_n$ is generated. This assumption does not require $\nu_n$ to exactly minimise the KL divergence and can be verified by leveraging convergence guarantees for variational inference algorithms. See for instance \cite{domke2023, kim2023} for gradient descent algorithms over location-scale families, \cite{lambert2022} for Wasserstein gradient descent algorithms over Gaussian densities, or \cite{wu2024} for natural gradient algorithms over exponential families.

\section{A logarithmic-moment approach for verifying Assumption (A1)}

\label{sec:log-moments}

This section presents sufficient conditions for Assumption (A1) to hold, that is, for $m_{n}=\inf_{\nu\in\mathcal{F}_{n}}\mathrm{KL}(\nu\Vert\mu_{n})$
to remain bounded over a family of approximating probability measures
$\mathcal{F}_{n}$. We do this by constructing a sequence $\nu_{n}^{*}\in\mathcal{F}_{n}$
such that $\limsup_{n\to\infty}\mathrm{KL}(\nu_{n}^{*}\Vert\mu_{n})<\infty$. Proofs of the results in this section and Section~\ref{sec:GLM} are collected in Appendix~\ref{sec:proofs}.

In the Euclidean setting, let $\lambda$ denote the Lebesgue measure on
$\mathbb{R}^{p}$ and let $\nu_{n}\ll\mu_{n}$, with $\mu_n$ equivalent to $\lambda$.
Then
\begin{equation}
    \mathrm{KL}(\nu_{n}\Vert\mu_{n})=\int\log\left(\frac{{\rm d}\nu_{n}}{{\rm d}\lambda}\right){\rm d}\nu_{n}-\int\log\left(\frac{{\rm d}\mu_{n}}{{\rm d}\lambda}\right){\rm d}\nu_{n}. \label{eq:KL_expansion}
\end{equation}
To bound $\mathrm{KL}(\nu_{n}\Vert\mu_{n})$, we must then
control log-moments of the densities of $\nu_{n}$ and $\mu_{n}$
under $\nu_{n}$.
Direct application of~\eqref{eq:KL_expansion} requires control of expectations under varying measures, which is often technically delicate.
As a simplifying device we use measurable bijections to convert this to the problem of bounding the log-moments of the densities of $\nu_{n}$ and $\mu_{n}$ under a single measure~$\tilde\nu$.

\subsection{Setup}

Let each $\mu_{n}$ be a Borel probability measure on $\mathbb{R}^{p}$,
and fix Borel measurable bijections $\gamma_{n}:\mathbb{R}^{p}\to\mathbb{R}^{p}$.
Define the push-forward measures $\tilde\mu_n=\mu_{n}\circ\gamma_{n}^{-1}$.
We write $\tilde p_{n}$ for the Lebesgue density of $\tilde\mu_n$.

\begin{prop}
\label{prop:log-mom-oneQ} Let $\tilde\nu$ be a Borel probability measure on $\mathbb{R}^{p}$ with Lebesgue density $\tilde q$. Assume that, for some $n_{0}$ and for all $n\ge n_{0}$, the push-forward
$\nu_{n}^{*}=\tilde{\nu}\circ\gamma_{n}$ belongs to $\mathcal{F}_{n}$ and
satisfies $\nu_{n}^{*}\ll\mu_{n}$. Then, for
all $n\ge n_{0}$,
\[
\mathrm{KL}(\nu_{n}^{*}\Vert\mu_{n})=\mathrm{KL}(\tilde{\nu}\Vert\tilde{\mu}_{n})=\mathrm{E}_{\tilde{\nu}}\left[\log\left(\frac{\tilde{q}}{\tilde{p}_{n}}\right)\right].
\]
In particular, if $\mathrm{E}_{\tilde{\nu}}|\log \tilde{q}|<\infty$
and $\limsup_{n\to\infty}\mathrm{E}_{\tilde{\nu}}(-\log \tilde{p}_{n})_+<\infty$,
then (A1) holds.
\end{prop}

\begin{rem}
\label{rem:renyi-feasibility-location-scale} The feasibility condition
$\nu_{n}^{*}\in\mathcal{F}_{n}$ for all large $n$ is usually easy
to verify for common variational classes. A common choice of bijection
is the affine recentering and rescaling $\gamma_{n}(\theta)=t_{n}(\theta-\theta_{n})$,
with centring points $\theta_{n}\in\mathbb{R}^{p}$ and positive scales
$t_{n}\to\infty$. Then, for any probability measure $\tilde\nu$ on $\mathbb{R}^{p}$
and random variable $Z\sim \tilde\nu$, $\tilde\nu\circ\gamma_{n}$ is the measure of the law $\mathcal{L}(\theta_{n}+t_{n}^{-1}Z)$.
Suppose there is a class $\mathcal{V}$ of measures on $\mathbb{R}^{p}$
(e.g.\ a location-scale or elliptical family) such that, for each
$n$,
\[
\mathcal{F}_{n}\supseteq\left\{ \mathcal{L}\!\left(\theta_{n}+t_{n}^{-1}Z\right):Z\sim \tilde\nu,\ \tilde\nu\in\mathcal{V}\right\} .
\]
Then, for any fixed $\tilde\nu\in\mathcal{V}$, the reference measure $\nu_{n}^{*}=\tilde\nu\circ\gamma_{n}$
belongs to $\mathcal{F}_{n}$ for all $n$. Moreover, if $\mu_{n}$
has Lebesgue density $p_{n}$, then the corresponding push-forward
$\tilde\mu_n=\mu_{n}\circ\gamma_{n}^{-1}$ has density
\[
\tilde p_{n}(z)=t_{n}^{-p}p_{n}\!\left(\theta_{n}+t_{n}^{-1}z\right).
\]
\end{rem}

\subsection{A simple sufficient envelope condition for positive-part log-moments}

A convenient way to verify the positive-part log-moment bound $\limsup_{n}\mathrm{E}_{\tilde\nu}(-\log \tilde p_{n})_+<\infty$
is to exhibit a pointwise lower envelope for $\tilde p_{n}$. Since $x\mapsto-\log x$
is decreasing, any bound $\tilde p_{n} \geq \exp (-h)$ implies
\[
(-\log \tilde p_{n})_+\le h_+.
\]
In particular, it is enough that
$\tilde p_{n}(z)\ge\exp\left(-h(z)\right)$ for all $n\ge n_{0}$ and some measurable $h$ with $\mathrm{E}_{\tilde\nu}h_+<\infty$.
This shows that the positive-part log-moment condition only prevents $\tilde p_{n}$ from
becoming too small on sets carrying non-negligible $\tilde\nu$-mass. The
following concrete application of this principle is used in the sequel.

\begin{prop}
\label{prop:LAN-like} If there exist $c,C>0$ and $n_{0}\in\mathbb N$ such that $\tilde p_{n}(z)\ge c\exp\left(-\frac{C}{2}\|z\|^{2}\right)$
for all $z\in\mathbb{R}^{p}$ and all $n\ge n_0$, and $\tilde\nu$ has finite second moment,
then $\limsup_{n\to\infty}\mathrm{E}_{\tilde\nu}(-\log \tilde p_{n})_+<\infty$.
\end{prop}

The next result isolates a mechanism that yields the required lower envelope.

\begin{prop} \label{prop:quadratic-envelope} Let $\tilde\mu_{n}$ be Borel
probability measures on $\mathbb{R}^{p}$ with strictly positive Lebesgue
densities $\tilde p_{n}$. Write $g_{n}=\log \tilde p_{n}$. Assume that there exist
constants $L,M,\alpha,r>0$ and $n_{0}\in\mathbb{N}$ such that, for
every $n\ge n_{0}$, $\mathrm{(B1)}$  $\nabla g_{n}$ is $L$-Lipschitz on $\mathbb{R}^{p}$, $\mathrm{(B2)}$  $\|\nabla g_{n}(0)\|\le M$, $\mathrm{(B3)}$ $\tilde\mu_{n}(B(0,r))\ge\alpha$.
Then there exist constants $c,C>0$, depending only on $L,M,\alpha,r$
and $p$, such that
\[
\tilde p_{n}(z)\ge c\exp\left(-\frac{C}{2}\|z\|^{2}\right)\qquad\text{for all }z\in\mathbb{R}^{p}\text{ and all }n\ge n_{0}.
\]
\end{prop}

A sufficient condition for $\mathrm{(B3)}$ is the existence of a
weak limit $\mu_{\infty}$ of $(\tilde\mu_n)$ with $\mu_{\infty}(B(0,r))>0$
for some $r>0$. Indeed, the Portmanteau theorem then yields
\[
\liminf_{n\to\infty}\tilde\mu_{n}(B(0,r))\ge\mu_{\infty}(B(0,r))>0,
\]
so any $0<\alpha<\mu_{\infty}(B(0,r))$ satisfies assumption $\mathrm{(B3)}$
for all sufficiently large $n$.

\subsection{Satisfying the envelope condition in the case of generalised posteriors}

\label{sec:stochasticVerification}

Let $(\Omega,\mathfrak{A},\Pr)$ be the data-generating probability space, with expectation $\mathrm{E}$, and let $(X_{i})_{i\ge1}$
be i.i.d.\ $\mathbb{X}$-valued random elements with common law $\mathrm{P}_{0}$.
This probability notation is used for the random data, while $\mu$ and $\nu$ continue to denote Borel probability measures on parameter spaces. Given a loss $\ell:\mathbb{R}^{p}\times\mathbb{X}\to\mathbb{R}$,
a prior density $\pi$ on $\mathbb{R}^{p}$, and temperatures $\eta_{n}>0$,
consider the generalised posterior \citep{bissiri2016general}
\begin{align}
\mu_{n}(\mathrm{d}\theta)\propto\exp\left\{ -\eta_{n}\sum_{i=1}^{n}\ell(\theta,X_{i})\right\} \pi(\theta)\,\mathrm{d}\theta,\label{eq:Gibbs-posterior}
\end{align}
with random Lebesgue density $p_{n}$. Let $\hat{\theta}_{n}:\Omega\to\mathbb{R}^{p}$
be a measurable centring sequence. For each $\omega\in\Omega$, define
\[
\gamma_{n}^{\omega}(\theta)=\sqrt{n}\,(\theta-\hat{\theta}_{n}(\omega)),\qquad\tilde\mu_n(\omega)=\mu_{n}(\omega)\circ(\gamma_{n}^{\omega})^{-1},
\]
and let $\tilde p_{n}(\omega,\cdot)$ denote the Lebesgue density of $\tilde\mu_n(\omega)$,
so that
\[
\tilde p_{n}(\omega,z)=n^{-p/2}p_{n}\!\left(\omega,\hat{\theta}_{n}(\omega)+n^{-1/2}z\right).
\]
Below, $B(x,r)$ and $\overline{B}(x,r)$ are understood with respect to the Euclidean metric on $\mathbb{R}^{p}$.
The next result verifies the hypotheses
of Proposition~\ref{prop:quadratic-envelope} pathwise.

\begin{thm} \label{thm:gp-envelope} Assume:
\begin{enumerate}
\item[$\mathrm{(C1)}$] For every $x\in\mathbb{X}$, the map $\theta\mapsto\ell(\theta,x)$
is twice continuously differentiable, and there exists a measurable
$H:\mathbb{X}\to[0,\infty)$ such that
\[
\sup_{\theta\in\mathbb{R}^{p}}\|\nabla_{\theta}^{2}\ell(\theta,x)\|\le H(x),\qquad\mathrm{E}[H(X_{1})]<\infty.
\]
\item[$\mathrm{(C2)}$] $\pi$ is strictly positive and twice continuously differentiable
on $\mathbb{R}^{p}$, and
\[
M_{\pi}=\sup_{\theta\in\mathbb{R}^{p}}\|\nabla_{\theta}^{2}\log\pi(\theta)\|<\infty.
\]
\item[$\mathrm{(C3)}$] $\bar{\eta}=\sup_{n\ge1}\eta_{n}<\infty$.
\item[$\mathrm{(C4)}$] There exist almost sure events $\Omega_{\mathrm{loss}},\Omega_{\nabla}\in\mathfrak{A}$
such that, for every $\omega\in\Omega_{\mathrm{loss}}$ and every
$n\ge1$, $\hat{\theta}_{n}(\omega)$ is a global minimiser of $\theta\mapsto\sum_{i=1}^{n}\ell(\theta,X_{i}(\omega))$,
and, for every $\omega\in\Omega_{\nabla}$,
\[
n^{-1/2}\left\|\nabla_{\theta}\log\pi\left(\hat{\theta}_{n}(\omega)\right)\right\|\to0.
\]
\item[$\mathrm{(C5)}$] There exist a Borel probability measure $\mu_{\infty}$ on $\mathbb{R}^{p}$,
a radius $r>0$, and an event $\Omega_{\mathrm{w}}\in\mathfrak{A}$
with $\Pr(\Omega_{\mathrm{w}})=1$ such that, for every $\omega\in\Omega_{\mathrm{w}}$,
\[
\tilde\mu_n(\omega)\rightsquigarrow\mu_{\infty}\qquad\text{and}\qquad\mu_{\infty}(B(0,r))>0.
\]
\end{enumerate}
Then there exist constants $c,C>0$ and an event $\Omega_{0}\in\mathfrak{A}$
with $\Pr(\Omega_{0})=1$ such that, for every $\omega\in\Omega_{0}$,
there exists $n_{0}(\omega)\in\mathbb{N}$ for which
\begin{equation}
\tilde p_{n}(\omega,z)\ge c\exp\left(-\frac{C}{2}\|z\|^{2}\right)\qquad\text{for all }z\in\mathbb{R}^{p}\text{ and all }n\ge n_{0}(\omega).\label{eq:Gen_Bayes_Lower_bound}
\end{equation}
\end{thm}

\section{Application to generalised linear models in generalised Bayes}

\label{sec:GLM}

We now specialise the generalised posterior (\ref{eq:Gibbs-posterior})
to canonical generalised linear models (GLMs) for which the loss is
defined on all of $\mathbb{R}^{p}$. Let $X_{i}=(Y_{i},W_{i})$ be
i.i.d.\ with law $\mathrm{P}_{0}$, where $Y_{i}:\Omega\to\mathbb{R}$
and $W_{i}:\Omega\to\mathbb{R}^{p}$. Assume that, conditional on
$W_{i}=w$, the response $Y_{i}$ has a canonical one-parameter exponential-family
conditional density with natural parameter $\theta^{\top}w$, namely
\[
p_{\theta}(y\mid w)=\exp\{y\theta^{\top}w-b(\theta^{\top}w)\},\qquad\theta\in\mathbb{R}^{p},
\]
where $b$ is three times continuously differentiable and $b''(u)>0$ for all $u\in\mathbb R$.
Up to an additive term depending only on $y$, the associated loss
is
\[
\ell(\theta,(y,w))=b(\theta^{\top}w)-y\theta^{\top}w.
\]
Taking $\tilde\nu=\mathrm{N}(0,\mathbf{I}_{p})$, the next theorem verifies
the hypotheses needed to apply Theorem~\ref{thm:gp-envelope}, Proposition~\ref{prop:LAN-like},
and the positive-part logarithmic-moment step in Proposition~\ref{prop:log-mom-oneQ}
pathwise in this GLM setting.

Assume that the prior density $\pi$ and the measurable centring sequence
$(\hat{\theta}_{n})$ satisfy {(C2)} and {(C4)} of Theorem~\ref{thm:gp-envelope},
that $\eta_{n}\to\eta_{*}\in(0,\infty)$, and that
\[
H(X_{1})=\sup_{\theta\in\mathbb{R}^{p}}\|\nabla_{\theta}^{2}\ell(\theta,(Y_{1},W_{1}))\|\le\sup_{\theta\in\mathbb{R}^{p}}b''(\theta^{\top}W_{1})\|W_{1}\|^{2}\qquad\text{satisfies}\qquad\mathrm{E}[H(X_{1})]<\infty.
\]
Assume moreover that there exist $\theta_{0}\in\mathbb{R}^{p}$ and
$r_{0}>0$ such that
\[
\mathrm{E}\!\left[W_{1}\{b'(\theta_{0}^{\top}W_{1})-Y_{1}\}\right]=0,\qquad\mathrm{E}[\|W_{1}\|\,|Y_{1}|]<\infty,\qquad\mathrm{E}|b(\theta^{\top}W_{1})|<\infty\quad\text{for all }\theta\in\mathbb{R}^{p},
\]
\[
\mathrm{E}[W_{1}W_{1}^{\top}]\text{ exists and is positive definite},\qquad\mathrm{E}\left[\sup_{\theta\in\overline{B}(\theta_{0},r_{0})}|b'''(\theta^{\top}W_{1})|\,\|W_{1}\|^{3}\right]<\infty.
\]
\begin{thm}
    \label{thm:GLM} Under the preceding assumptions, the
hypotheses of Theorem~\ref{thm:gp-envelope} hold. Consequently,
there exist constants $c,C>0$ and an event $\Omega_{0}\in\mathfrak{A}$
with $\Pr(\Omega_{0})=1$ such that, for every $\omega\in\Omega_{0}$,
there exists $n_{0}(\omega)\in\mathbb{N}$ for which \eqref{eq:Gen_Bayes_Lower_bound}
holds.
\end{thm}

Gaussian location-scale variational classes satisfy the feasibility
condition $\nu_{n}^{*}=\tilde\nu\circ\gamma_{n}\in\mathcal{F}_{n}$ for all
large $n$. In particular, Theorem~\ref{thm:GLM} applies to
the Gaussian variational GLM approximations considered by \citet{JaakkolaJordan2000,KnowlesMinka2011}.

\newpage

\appendix

\section{Proofs}
\label{sec:proofs}

\begin{proof}[of Proposition~\ref{prop:log-mom-oneQ}]
By \citet[Cor. 2.18]{polyanskiy_wu_2025_information_theory}, $\mathrm{KL}(\nu_{n}^{*}\Vert\mu_{n})=\mathrm{KL}(\tilde\nu\Vert\tilde\mu_n)$,
and bijectivity of $\gamma_{n}$ implies $\tilde\nu\ll\tilde\mu_n$. Writing
$\tilde p_{n}={\rm d}\tilde\mu_n/{\rm d}\lambda$ and $\tilde q={\rm d}\tilde\nu/{\rm d}\lambda$,
the chain rule for Radon--Nikodym derivatives gives that on $\{\tilde p_{n}\neq 0\}$, ${\rm d}\tilde\nu/{\rm d}\tilde\mu_n=\tilde q/\tilde p_{n}$. As $\{\tilde p_{n}=0\}$ is a $\tilde\nu$ null set,
\[
\mathrm{KL}(\tilde\nu\Vert\tilde\mu_n)=\int\log\left(\frac{\tilde q}{\tilde p_{n}}\right){\rm d}\tilde\nu=\mathrm{E}_{\tilde\nu}\left[\log\left(\frac{\tilde q}{\tilde p_{n}}\right)\right].
\]
If the two sufficient log-moment bounds hold, then
\[
\log\left(\frac{\tilde q}{\tilde p_n}\right)\le |\log \tilde q|+(-\log \tilde p_n)_+,
\]
and therefore
\[\limsup_{n\to\infty} m_n \leq \limsup_{n\to\infty} \mathrm{KL}(\nu_{n}^{*}\Vert\mu_{n}) < \infty.\]
This establishes (A1).
\end{proof}

\begin{proof}[of Proposition~\ref{prop:LAN-like}]
The bound on $\tilde p_{n}$ implies $-\log \tilde p_{n}(\cdot)\le-\log c+\frac{C}{2}\|\cdot\|^{2}$,
hence uniformly over $n$,
\[
\mathrm{E}_{\tilde\nu}(-\log \tilde p_{n})_+\le |\log c|+\frac{C}{2}\int\| z\|^{2}~\tilde{\nu}(\mathrm{d}z)<\infty.
\]

\end{proof}

\begin{proof}[of Proposition~\ref{prop:quadratic-envelope}] Fix $n\ge n_{0}$.
Since $\nabla g_{n}$ is $L$-Lipschitz, for every $z\in\mathbb{R}^{p}$,
\[
g_{n}(z)=g_{n}(0)+\int_{0}^{1}\langle\nabla g_{n}(tz),z\rangle\,\mathrm{d}t.
\]
Hence
\[
\begin{aligned}g_{n}(z) & =g_{n}(0)+\langle\nabla g_{n}(0),z\rangle+\int_{0}^{1}\langle\nabla g_{n}(tz)-\nabla g_{n}(0),z\rangle\,\mathrm{d}t\\
 & \ge g_{n}(0) - \|\nabla g_n(0)\|\, \|z\| -\int_{0}^{1}\|\nabla g_{n}(tz)-\nabla g_{n}(0)\|\,\|z\|\,\mathrm{d}t\\
 & \ge g_{n}(0)-M\|z\|-\int_{0}^{1}Lt\,\|z\|^{2}\,\mathrm{d}t\\
 & \ge g_{n}(0)-M\|z\|-\frac{L}{2}\|z\|^{2}.
\end{aligned}
\]
Therefore
\begin{equation}
\tilde p_{n}(z)\ge \tilde p_{n}(0)\exp\left(-M\|z\|-\frac{L}{2}\|z\|^{2}\right).\label{eq:quadratic-envelope-lower-step}
\end{equation}

Applying the same argument but instead using that $\langle\nabla g_{n}(x),z\rangle \leq \|\nabla g_n(x)\|\,\|z\|$, for every $x\in\mathbb{R}^{p}$,
\[
g_{n}(x)\le g_{n}(0)+M\|x\|+\frac{L}{2}\|x\|^{2}.
\]
Hence, for every $x\in B(0,r)$,
\[
\tilde p_{n}(x)\le \tilde p_{n}(0)\exp\left(Mr+\frac{L}{2}r^{2}\right).
\]
Integrating over $B(0,r)$ and using assumption~(B3)
gives
\[
\alpha\le\tilde\mu_{n}(B(0,r))=\int_{B(0,r)}\tilde p_{n}(x)\,\mathrm{d}x\le\lambda\left(B(0,r)\right)\exp\left(Mr+\frac{L}{2}r^{2}\right)\tilde p_{n}(0),
\]
so that
\[
\tilde p_{n}(0)\ge\frac{\alpha}{\lambda\left(B(0,r)\right)}\exp\left(-Mr-\frac{L}{2}r^{2}\right)=:m_{0}>0.
\]
Substituting this bound into \eqref{eq:quadratic-envelope-lower-step}
yields
\[
\tilde p_{n}(z)\ge m_{0}\exp\left(-M\|z\|-\frac{L}{2}\|z\|^{2}\right).
\]
Finally, by Young's inequality,
\[
M\|z\|\le\frac{M^{2}}{2}+\frac{\|z\|^{2}}{2},
\]
so
\[
\tilde p_{n}(z)\ge m_{0}\,e^{-M^{2}/2}\exp\left(-\frac{L+1}{2}\|z\|^{2}\right).
\]
Thus the conclusion holds with
\[
c=\frac{\alpha}{\lambda\left(B(0,r)\right)}\exp\left(-Mr-\frac{L}{2}r^{2}-\frac{M^{2}}{2}\right),\qquad C=L+1.
\]
\end{proof}

\begin{proof}[of Theorem~\ref{thm:gp-envelope}] By the strong
law of large numbers,
\[
\Omega_{H}=\left\{\omega\in\Omega: \frac{1}{n}\sum_{i=1}^{n}H(X_{i}(\omega))\to\mathrm{E}[H(X_{1})]\right\}
\]
has probability one. Let
\[
\Omega_{0}=\Omega_{H}\cap\Omega_{\mathrm{loss}}\cap\Omega_{\nabla}\cap\Omega_{\mathrm{w}},
\]
so that $\Pr(\Omega_{0})=1$. Fix $\omega\in\Omega_{0}$ for the remainder
of the proof and suppress the $\omega$-dependence when no confusion
arises.

For each $n$, define
\[
g_{n}(z)=\log \tilde p_{n}(z),\qquad z\in\mathbb{R}^{p}.
\]
Since
\[
g_{n}(z)=-\eta_{n}\sum_{i=1}^{n}\ell\!\left(\hat{\theta}_{n}+n^{-1/2}z,X_{i}\right)+\log\pi\!\left(\hat{\theta}_{n}+n^{-1/2}z\right)-\log\mathcal{Z}_{n}-\frac{p}{2}\log n,
\]
for a normalising constant $\mathcal{Z}_{n}$ independent of $z$,
assumptions {(C1)} and {(C2)} imply that $g_{n}$ is twice
continuously differentiable and
\[
\nabla_{z}^{2}g_{n}(z)=\frac{1}{n}\left\{ -\eta_{n}\sum_{i=1}^{n}\nabla_{\theta}^{2}\ell\!\left(\hat{\theta}_{n}+n^{-1/2}z,X_{i}\right)+\nabla_{\theta}^{2}\log\pi\!\left(\hat{\theta}_{n}+n^{-1/2}z\right)\right\} .
\]
Hence
\[
\sup_{z\in\mathbb{R}^{p}}\|\nabla_{z}^{2}g_{n}(z)\|\le\eta_{n}\frac{1}{n}\sum_{i=1}^{n}H(X_{i})+\frac{M_{\pi}}{n}.
\]
Since $\omega\in\Omega_{H}$, there exists $n_{1}(\omega)\in\mathbb{N}$
such that, with $L=\bar{\eta}\left(\mathrm{E}[H(X_{1})]+1\right)+M_{\pi}$,
\[
\sup_{z\in\mathbb{R}^{p}}\|\nabla_{z}^{2}g_{n}(z)\|\le L\qquad\text{for all }n\ge n_{1}(\omega).
\]
In particular, $\nabla g_{n}$ is $L$-Lipschitz on $\mathbb{R}^{p}$
for all $n\ge n_{1}(\omega)$.

Next,
\[
\nabla_{z}g_{n}(0)=n^{-1/2}\left\{ -\eta_{n}\sum_{i=1}^{n}\nabla_{\theta}\ell(\hat{\theta}_{n},X_{i})+\nabla_{\theta}\log\pi(\hat{\theta}_{n})\right\} .
\]
Because $\hat{\theta}_{n}$ is a global minimiser of $\theta\mapsto\sum_{i=1}^{n}\ell(\theta,X_{i})$
and this criterion is continuously differentiable on $\mathbb{R}^{p}$,
we have
\[
\sum_{i=1}^{n}\nabla_{\theta}\ell(\hat{\theta}_{n},X_{i})=0.
\]
Therefore
\[
\nabla_{z}g_{n}(0)=n^{-1/2}\nabla_{\theta}\log\pi(\hat{\theta}_{n}).
\]
Since $\omega\in\Omega_{\nabla}$, there exists $n_{2}(\omega)\in\mathbb{N}$
such that
\[
\|\nabla_{z}g_{n}(0)\|\le1\qquad\text{for all }n\ge n_{2}(\omega).
\]

Since $\omega\in\Omega_{\mathrm{w}}$, we have
\[
\tilde\mu_n\rightsquigarrow\mu_{\infty}\qquad\text{and}\qquad\mu_{\infty}(B(0,r))>0.
\]
The Portmanteau theorem gives
\[
\liminf_{n\to\infty}\tilde\mu_n(B(0,r))\ge\mu_{\infty}(B(0,r))>0.
\]
Hence, with
\[
\alpha=\frac{1}{2}\mu_{\infty}(B(0,r))>0,
\]
there exists $n_{3}(\omega)\in\mathbb{N}$ such that
\[
\tilde\mu_n(B(0,r))\ge\alpha\qquad\text{for all }n\ge n_{3}(\omega).
\]

Applying Proposition~\ref{prop:quadratic-envelope} to the deterministic
sequence $(\tilde\mu_n(\omega))_{n\ge1}$ with the constants $L$, $M=1$,
$\alpha$, and $r$, we obtain constants $c,C>0$, depending only
on $\bar{\eta},~\mathrm{E}[H(X_{1})],~M_{\pi},~r,~\mu_{\infty}(B(0,r)),~p$,
and an index
\[
n_{0}(\omega)=\max\{n_{1}(\omega),n_{2}(\omega),n_{3}(\omega)\}
\]
such that
\[
\tilde p_{n}(\omega,z)\ge c\exp\left(-\frac{C}{2}\|z\|^{2}\right)\qquad\text{for all }z\in\mathbb{R}^{p}\text{ and all }n\ge n_{0}(\omega).
\]
\end{proof}

\begin{proof}[of Theorem~\ref{thm:GLM}] We verify {(C1)},
{(C3)}, and {(C5)} of Theorem~\ref{thm:gp-envelope}; assumptions
{(C2)} and {(C4)} are part of the statement. Write
\[
B_0=B(\theta_{0},r_{0}).
\]

For $(y,w)\in\mathbb{R}\times\mathbb{R}^{p}$,
\[
\nabla_{\theta}\ell(\theta,(y,w))=\{b'(\theta^{\top}w)-y\}w,\qquad\nabla_{\theta}^{2}\ell(\theta,(y,w))=b''(\theta^{\top}w)\,ww^{\top}.
\]
Hence
\[
\|\nabla_{\theta}^{2}\ell(\theta,(Y_{1},W_{1}))\|=\|b''(\theta^{\top}W_{1})W_{1}W_{1}^{\top}\|\le H(X_{1}),
\]
since $\|W_{1}W_{1}^{\top}\|=\|W_{1}\|^{2}$. Therefore {(C1)}
holds. Since $\eta_{n}\to\eta_{*}\in(0,\infty)$, we also have $\bar{\eta}=\sup_{n\ge1}\eta_{n}<\infty$,
so {(C3)} holds.

To verify {(C5)}, set
\[
g_{n}(\theta)=\frac{1}{n}\sum_{i=1}^{n}\ell(\theta,X_{i})=\frac{1}{n}\sum_{i=1}^{n}\{b(\theta^{\top}W_{i})-Y_{i}\theta^{\top}W_{i}\},\qquad g(\theta)=\mathrm{E}[\ell(\theta,X_{1})],
\]
and define
\[
\tilde{g}_{n}(\theta)=\eta_{n}g_{n}(\theta),\qquad\tilde{g}(\theta)=\eta_{*}g(\theta).
\]
Since
\[
-\log p_{\theta}(Y_{i}\mid W_{i})=b(\theta^{\top}W_{i})-Y_{i}\theta^{\top}W_{i},
\]
up to an additive term independent of $\theta$, the criterion $g_{n}$
is exactly the empirical criterion appearing in \citet[Thm.~13]{miller2021generalized_posteriors}
for canonical GLMs. The assumptions in the statement verify
the hypotheses of that theorem. The only points requiring verification
are identifiability and the third-derivative bound.

First, if $a\in\mathbb{R}^{p}$ satisfies $a^{\top}W_{1}=0$ almost
surely, then
\[
\mathrm{E}[(a^{\top}W_{1})^{2}]=a^{\top}\mathrm{E}[W_{1}W_{1}^{\top}]a=0.
\]
Since $\mathrm{E}[W_{1}W_{1}^{\top}]$ is positive definite, it follows
that $a=0$. Thus the identifiability condition in \citet[Thm.~13]{miller2021generalized_posteriors}
holds. Second, for every $j,k,\ell\in\{1,\ldots,p\}$,
\[
|W_{1j}W_{1k}W_{1\ell}|\le\|W_{1}\|^{3},
\]
and hence
\[
\mathrm{E}\!\left[\sup_{\theta\in\overline{B}(\theta_{0},r_{0})}|b'''(\theta^{\top}W_{1})|\,|W_{1j}W_{1k}W_{1\ell}|\right]\le\mathrm{E}\!\left[\sup_{\theta\in\overline{B}(\theta_{0},r_{0})}|b'''(\theta^{\top}W_{1})|\,\|W_{1}\|^{3}\right]<\infty.
\]

Therefore \citet[Thm.~13]{miller2021generalized_posteriors} implies
that, on an event $\Omega_{\mathrm{M}}\in\mathfrak{A}$ with $\Pr(\Omega_{\mathrm{M}})=1$,
the sequence $(g_{n}(\omega,\cdot))_{n\ge1}$ satisfies the hypotheses
of case~(2) of \citet[Thm.~5]{miller2021generalized_posteriors}
for every $\omega\in\Omega_{\mathrm{M}}$. Since $\eta_{n}\to\eta_{*}\in(0,\infty)$,
it follows that, for every $\omega\in\Omega_{\mathrm{M}}$, the sequence
$(\tilde{g}_{n}(\omega,\cdot))_{n\ge1}$ also satisfies the hypotheses
of case~(2) of \citet[Thm.~5]{miller2021generalized_posteriors},
with limit $\tilde{g}$. In particular, by \citet[Thm.~7]{miller2021generalized_posteriors},
for every $\omega\in\Omega_{\mathrm{M}}$,
\[
\tilde{g}_{n}(\omega,\cdot)\to\tilde{g}\qquad\text{and}\qquad\nabla_{\theta}^{2}\tilde{g}_{n}(\omega,\cdot)\to\nabla_{\theta}^{2}\tilde{g}
\]
uniformly on $B_0$.

Let
\[
H_{3}(X_{1})=\sup_{\theta\in\overline{B}(\theta_{0},r_{0})}|b'''(\theta^{\top}W_{1})|\,\|W_{1}\|^{3}.
\]
By assumption, $\mathrm{E}[H_{3}(X_{1})]<\infty$. Hence the strong
law of large numbers yields an event $\Omega_{3}\in\mathfrak{A}$
with $\Pr(\Omega_{3})=1$ such that
\[
\frac{1}{n}\sum_{i=1}^{n}H_{3}(X_{i})\to\mathrm{E}[H_{3}(X_{1})]\qquad\text{on }\Omega_{3}.
\]

Fix $\omega\in\Omega_{\mathrm{M}}\cap\Omega_{3}\cap\Omega_{\mathrm{loss}}$,
where $\Omega_{\mathrm{loss}}$ is the event in the first part of
{(C4)}, and suppress the $\omega$-dependence for the remainder
of the proof. Since $\eta_{n}>0$, the minimisers of
\[
\theta\mapsto\sum_{i=1}^{n}\ell(\theta,X_{i}),\qquad g_{n}(\cdot),\qquad\text{and}\qquad\tilde{g}_{n}(\cdot)
\]
coincide. Hence $\hat{\theta}_{n}$ is a global minimiser of $\tilde{g}_{n}$
for every $n$. Because $\mathbb{R}^{p}$ is open and $\tilde{g}_{n}$
is continuously differentiable, we have
\[
\nabla_{\theta}\tilde{g}_{n}(\hat{\theta}_{n})=0\qquad\text{for every }n.
\]

We next show that $\hat{\theta}_{n}\to\theta_{0}$. Since case~(2)
of \citet[Thm.~5]{miller2021generalized_posteriors} implies case~(1),
there exists a compact set $\mathbb{K}\subset B_0$ with
$\theta_{0}\in\mathrm{int}(\mathbb{K})$ such that
\[
\tilde{g}(\theta)>\tilde{g}(\theta_{0})\qquad\text{for all }\theta\in\mathbb{K}\setminus\{\theta_{0}\},
\]
and
\[
\liminf_{n\to\infty}\inf_{\theta\in\mathbb{R}^{p}\setminus\mathbb{K}}\tilde{g}_{n}(\theta)>\tilde{g}(\theta_{0}).
\]
Choose $\rho>0$ such that
\[
\overline{B}(\theta_{0},\rho)\subset\mathrm{int}(\mathbb{K}).
\]
Since $\tilde{g}$ is continuous on $B_0$ and strictly
larger than $\tilde{g}(\theta_{0})$ on $\mathbb{K}\setminus\{\theta_{0}\}$,
there exists $\delta>0$ such that
\[
\inf_{\theta\in\mathbb{K}\setminus\overline{B}(\theta_{0},\rho)}\tilde{g}(\theta)\ge\tilde{g}(\theta_{0})+4\delta.
\]
By uniform convergence of $\tilde{g}_{n}$ to $\tilde{g}$ on $B_0$,
by convergence of $\tilde{g}_{n}(\theta_{0})$ to $\tilde{g}(\theta_{0})$,
and by the displayed liminf bound outside $\mathbb{K}$, there exists
$n_{1}\in\mathbb{N}$ such that, for all $n\ge n_{1}$,
\[
\tilde{g}_{n}(\theta_{0})\le\tilde{g}(\theta_{0})+\delta,
\]
\[
\inf_{\theta\in\mathbb{K}\setminus\overline{B}(\theta_{0},\rho)}\tilde{g}_{n}(\theta)\ge\tilde{g}(\theta_{0})+3\delta,\qquad\inf_{\theta\in\mathbb{R}^{p}\setminus\mathbb{K}}\tilde{g}_{n}(\theta)\ge\tilde{g}(\theta_{0})+3\delta.
\]
Hence
\[
\inf_{\theta\in\mathbb{R}^{p}\setminus\overline{B}(\theta_{0},\rho)}\tilde{g}_{n}(\theta)\ge\tilde{g}(\theta_{0})+3\delta>\tilde{g}_{n}(\theta_{0}).
\]
Since $\hat{\theta}_{n}$ is a global minimiser of $\tilde{g}_{n}$,
it follows that
\[
\hat{\theta}_{n}\in\overline{B}(\theta_{0},\rho)\qquad\text{for all }n\ge n_{1}.
\]
Because $\rho>0$ may be chosen arbitrarily small, we conclude that
\[
\hat{\theta}_{n}\to\theta_{0}.
\]

We now verify the hypotheses of \citet[Thm.~6]{miller2021generalized_posteriors}.
Condition~(2) holds because
\[
\nabla_{\theta}^{2}\tilde{g}_{n}(\theta_{0})\to\nabla_{\theta}^{2}\tilde{g}(\theta_{0}),
\]
and, for every $a\in\mathbb{R}^{p}\setminus\{0\}$,
\[
\begin{aligned}a^{\top}\nabla_{\theta}^{2}\tilde{g}(\theta_{0})a & =\eta_{*}\,a^{\top}\mathrm{E}\!\left[b''(\theta_{0}^{\top}W_{1})W_{1}W_{1}^{\top}\right]a\\
 & =\eta_{*}\,\mathrm{E}\!\left[b''(\theta_{0}^{\top}W_{1})(a^{\top}W_{1})^{2}\right]>0.
\end{aligned}
\]
Here we used that $\eta_{*}>0$, that $b''>0$ by assumption, and that $a^{\top}W_{1}$ is not almost surely zero by the
identifiability argument above. Hence $\nabla_{\theta}^{2}\tilde{g}(\theta_{0})$
is positive definite.

To verify Condition~(3), note that for every $j,k,\ell\in\{1,\ldots,p\}$,
\[
\sup_{\theta\in B_0}\left|\partial_{jk\ell}^{3}\tilde{g}_{n}(\theta)\right|\le\eta_{n}\frac{1}{n}\sum_{i=1}^{n}\sup_{u\in\overline{B}(\theta_{0},r_{0})}|b'''(u^{\top}W_{i})|\,|W_{ij}W_{ik}W_{i\ell}|\le\eta_{n}\frac{1}{n}\sum_{i=1}^{n}H_{3}(X_{i}).
\]
Since $\omega\in\Omega_{3}$ and $\eta_{n}\to\eta_{*}$, the right-hand
side is eventually bounded in $n$. Therefore the third derivatives
of $\tilde{g}_{n}$ are uniformly bounded on $B_0$, so
Condition~(3) of \citet[Thm.~6]{miller2021generalized_posteriors}
holds. Consequently, Assumption~(1) of \citet[Thm.~4]{miller2021generalized_posteriors}
is satisfied for $\tilde{g}_{n}$ with centring sequence $\hat{\theta}_{n}$.

To verify Assumption~(2) of \citet[Thm.~4]{miller2021generalized_posteriors},
fix $\varepsilon>0$. Since case~(2) of \citet[Thm.~5]{miller2021generalized_posteriors}
implies case~(1), and $\theta_{0}\in\mathrm{int}(\mathbb{K})$, there
exists $\delta>0$ such that
\[
\inf_{\theta\in\mathbb{K}\cap B(\theta_{0},\varepsilon/2)^{c}}\tilde{g}(\theta)\ge\tilde{g}(\theta_{0})+4\delta.
\]
By uniform convergence of $\tilde{g}_{n}$ to $\tilde{g}$ on $B_0$,
by convergence of $\tilde{g}_{n}(\theta_{0})$ to $\tilde{g}(\theta_{0})$,
and by the displayed liminf bound outside $\mathbb{K}$, there exists
$n_{2}\in\mathbb{N}$ such that, for all $n\ge n_{2}$,
\[
\tilde{g}_{n}(\theta_{0})\le\tilde{g}(\theta_{0})+\delta,
\]
\[
\inf_{\theta\in\mathbb{K}\cap B(\theta_{0},\varepsilon/2)^{c}}\tilde{g}_{n}(\theta)\ge\tilde{g}(\theta_{0})+3\delta,\qquad\inf_{\theta\in\mathbb{R}^{p}\setminus\mathbb{K}}\tilde{g}_{n}(\theta)\ge\tilde{g}(\theta_{0})+3\delta.
\]
Hence
\[
\inf_{\theta\in B(\theta_{0},\varepsilon/2)^{c}}\tilde{g}_{n}(\theta)\ge\tilde{g}(\theta_{0})+3\delta\qquad\text{for all }n\ge n_{2}.
\]
Since $\hat{\theta}_{n}\to\theta_{0}$, we have
\[
B(\theta_{0},\varepsilon/2)\subset B(\hat{\theta}_{n},\varepsilon)
\]
for all sufficiently large $n$. Therefore, for all sufficiently large
$n$,
\begin{align*}
\inf_{\theta\in B(\hat{\theta}_{n},\varepsilon)^{c}}\{\tilde{g}_{n}(\theta)-\tilde{g}_{n}(\hat{\theta}_{n})\} & \ge\inf_{\theta\in B(\theta_{0},\varepsilon/2)^{c}}\{\tilde{g}_{n}(\theta)-\tilde{g}_{n}(\hat{\theta}_{n})\}\\
 & \ge\inf_{\theta\in B(\theta_{0},\varepsilon/2)^{c}}\{\tilde{g}_{n}(\theta)-\tilde{g}_{n}(\theta_{0})\},
\end{align*}
since $\hat{\theta}_{n}$ minimises $\tilde{g}_{n}$. It follows from
the two preceding displays that
\[
\liminf_{n\to\infty}\inf_{\theta\in B(\hat{\theta}_{n},\varepsilon)^{c}}\{\tilde{g}_{n}(\theta)-\tilde{g}_{n}(\hat{\theta}_{n})\}>0.
\]
Thus Assumption~(2) of \citet[Thm.~4]{miller2021generalized_posteriors}
holds.

Since $\pi$ is strictly positive and twice continuously differentiable
by {(C2)}, the prior assumptions in \citet[Thm.~4]{miller2021generalized_posteriors}
are also satisfied. Finally,
\[
\mu_{n}(\mathrm{d}\theta)\propto\exp\{-n\tilde{g}_{n}(\theta)\}\pi(\theta)\,\mathrm{d}\theta.
\]
Hence \citet[Thm.~4]{miller2021generalized_posteriors} applies pathwise
and yields a nondegenerate centred Gaussian probability measure $\Phi$
on $\mathbb{R}^{p}$ such that
\[
\sup_{A\in\mathfrak{B}(\mathbb{R}^{p})}\left|\tilde\mu_n(A)-\Phi(A)\right|\to0.
\]
In particular, $\tilde\mu_n\rightsquigarrow\Phi$, and
\[
\Phi\left(B(0,r)\right)>0\qquad\text{for every }r>0.
\]
Thus {(C5)} of Theorem~\ref{thm:gp-envelope} holds on the almost
sure event
\[
\Omega_{\mathrm{w}}=\Omega_{\mathrm{M}}\cap\Omega_{3}\cap\Omega_{\mathrm{loss}}.
\]
Theorem~\ref{thm:gp-envelope} therefore yields the
conclusion.
\end{proof}

\end{document}